# 8 node solid brick element high order stiffness matrix template

Boning Zhang, Lan Nguyen
Department of Civil, Environmental, and Architectural Engineering
University of Colorado at Boulder

## 1 Introduction

More attention has been dedicated into finite element analysis (FEA) nowadays, especially to some analysis that can construct, process, and provide an efficient demonstration on an element performance. The intent of this paper is to develop a high order (HO) stiffness matrix template which is a parameterized algebraic representation of the element level stiffness equations that provide a continuum of consistent and stable finite element models of a given type and node/freedom configuration (Felippa, 2013) for an 8-node hexahedron (also known as brick) element. Template instances is produced by setting the value to free parameters furnish specific elements. Templates facilitate the unified implementation of finite element families, as well as the construction of custom elements (Felippa, 2013). In particular, the advantage of having a template is the symbolic power which can quickly solve and potentially lead to a comparison of how to construct an optimal element. The following discussion covers a brief background for an element stiffness matrix template.

In the early of 1970s, the concept of free formulation to templates was indirectly constructed. The element stiffness was to be derived directly from the consistency conditions and it was provided by the Individual Element Test of Bergan and Hanssen (Felippa, 2013) together with stability and accuracy considerations to determine algebraic redundancies if any. A modification in the Assumed Natural Strain (ANS) method of Park and Stanley was successfully tried for the HO stiffness. Needless to say, the method did not specifically address the HO stiffness for an eight-node solid element. The starting point for HO stiffness development could be found from the thin homogenous plate in plane stress shown in Felippa study in 2003, specifically in rectangular panel example, and could be also found in the stress element adopted from the study done by TCMT (short deviation for the name of four authors: Turner, Cloughm Martin, and Topp) back in 1956.

In this paper, the template will be developed from a Assumed Stress Method, which its formulation is based on the Hellinger-Reissner principle developed according to Kang's study in 1986. The element stiffness is decomposed into a basic part that takes care of consistency and mix-ability, and a HO element stiffness part that takes care of stability (also known as rank sufficient) and accuracy. In FE method, the HO stiffness is based on a displacement formulation, whereas the basis stiffness is method independent. To start, one should be familiar with the definition of a solid brick element. Solid brick element is three-dimensional (3D) finite elements that can model solid bodies and structures without any a priori geometric simplification (Felippa, 2013) . Figure 1 shows a simple sketch of a 3D hexahedron element.    The 8-node brick element contains twenty-four degrees of free-dom (DOF) collectively representing the linear displacements at each of the element nodes. In other words, a model would be provided with the config-uration of how the nodes connect together, and how its degree of freedoms link with in the chain of numbering system.

Figure 1. 3D hexahedron brick element

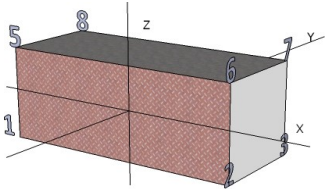

## 2 Assumed stress field

For a brick element subjected to axial and transverse loading, there are normal stresses, flexural stress, shear stress. One can image from the name *assumed stress method* that this method is based on an assumed stress field. From Kang's study, these assumed stress field for elements with parallel sides are as the following:

$$\sigma_{11} = \beta_1 + \beta_2\eta + \beta_3\mu + \beta_4\eta\mu$$
$$\sigma_{22} = \beta_5 + \beta_6\xi + \beta_7\mu + \beta_8\xi\mu$$
$$\sigma_{33} = \beta_9 + \beta_{10}\xi + \beta_{11}\eta + \beta_{12}\xi\eta$$
$$\sigma_{12} = \beta_{13} + \beta_{14}\mu$$
$$\sigma_{23} = \beta_{15} + \beta_{16}\xi$$
$$\sigma_{13} = \beta_{17} + \beta_{18}\eta$$

Notice that these stresses are proven to be satisfied the equilibrium equations with zero body forces according to Kang's study. Yet, a quick check could also be done by substituting these stress formulations into $\sigma_{ij,j} = 0$.

All these stresses are proven to be related to geometric of a given cubic size and a stress amplitude parameter matrix, the following shows the matrix written format:

$$\underbrace{\sigma}_{6\times 1} = \underbrace{N}_{6\times 18} \underbrace{\beta}_{18\times 1} \tag{1}$$

Note that in *equation (1)*, N matrix is not a shape function vector, rather it is a geometric matrix:

$$\mathbf{N} = \begin{bmatrix} 1 & \eta & \mu & \eta\mu & 0 & 0 & 0 & 0 & 0 & 0 & 0 & 0 & 0 & 0 & 0 & 0 \\ 0 & 0 & 0 & 0 & 1 & \xi & \mu & \mu\xi & 0 & 0 & 0 & 0 & 0 & 0 & 0 & 0 \\ 0 & 0 & 0 & 0 & 0 & 0 & 0 & 0 & 1 & \xi & \eta & \eta\xi & 0 & 0 & 0 & 0 \\ 0 & 0 & 0 & 0 & 0 & 0 & 0 & 0 & 0 & 0 & 0 & 1 & \mu & 0 & 0 & 0 \\ 0 & 0 & 0 & 0 & 0 & 0 & 0 & 0 & 0 & 0 & 0 & 0 & 1 & \xi & 0 & 0 \\ 0 & 0 & 0 & 0 & 0 & 0 & 0 & 0 & 0 & 0 & 0 & 0 & 0 & 0 & 1 & \eta \end{bmatrix}$$

The relationship between nodal force vector, **f**, and the parameter vector, **/3**, can be calculated based on *Cauchy traction formulation*: **CFn = t**.

As in *figure 1*, take *face 2376* as an example

$$\mathbf{t} = \begin{bmatrix} t_x \\ t_y \\ t_z \end{bmatrix} = \boldsymbol{\sigma}\mathbf{n} = \begin{bmatrix} \sigma_{11} & \sigma_{12} & \sigma_{13} \\ \sigma_{12} & \sigma_{22} & \sigma_{23} \\ \sigma_{13} & \sigma_{23} & \sigma_{33} \end{bmatrix} \begin{bmatrix} 1 \\ 0 \\ 0 \end{bmatrix} = \begin{bmatrix} \sigma_{11} \\ \sigma_{12} \\ \sigma_{13} \end{bmatrix}$$

Then from the assumed stress field:

$$t_x = \sigma_{11} = \beta_1 + \beta_2\eta + \beta_3\mu + \beta_4\eta\mu = \begin{bmatrix} \beta_1 & \beta_2 & \beta_3 & \beta_4 \end{bmatrix} \begin{bmatrix} 1 \\ \eta \\ \mu \\ \eta\mu \end{bmatrix}$$

Assume face displacement as:

$$u_x = u_{x2}N_2 + u_{x3}N_3 + u_{x7}N_7 + u_{x6}N_6 = \begin{bmatrix} N_2 & N_3 & N_7 & N_6 \end{bmatrix} \begin{bmatrix} u_{x2} \\ u_{x3} \\ u_{x7} \\ u_{x6} \end{bmatrix},$$

$$N_2 = \frac{1}{4}(1-\eta)(1-\mu), \quad N_3 = \frac{1}{4}(1+\eta)(1-\mu)$$
$$N_7 = \frac{1}{4}(1+\eta)(1+\mu), \quad N_6 = \frac{1}{4}(1-\eta)(1+\mu)$$

Then

$$\mathbf{f}_x^T \mathbf{u}_x = \int_{face2376} t_x \cdot u_x dA$$

$$= \int_{face2376} \begin{bmatrix} \beta_1 & \beta_2 & \beta_3 & \beta_4 \end{bmatrix} \begin{bmatrix} 1 \\ \eta \\ \mu \\ \eta\mu \end{bmatrix} \begin{bmatrix} N_2 & N_3 & N_7 & N_6 \end{bmatrix} \begin{bmatrix} u_{x2} \\ u_{x3} \\ u_{x7} \\ u_{x6} \end{bmatrix} dA$$

$$= \begin{bmatrix} \beta_1 & \beta_2 & \beta_3 & \beta_4 \end{bmatrix} \int_{-1}^{1}\int_{-1}^{1} \begin{bmatrix} 1 \\ \eta \\ \mu \\ \eta\mu \end{bmatrix} \begin{bmatrix} N_2 & N_3 & N_7 & N_6 \end{bmatrix} J_{face2376} d\eta d\mu \begin{bmatrix} u_{x2} \\ u_{x3} \\ u_{x7} \\ u_{x6} \end{bmatrix}$$

From which, the nodal force vector in x direction will be

$$\mathbf{f}_x^T = \begin{bmatrix} \beta_1 & \beta_2 & \beta_3 & \beta_4 \end{bmatrix} \int_{-1}^{1}\int_{-1}^{1} \begin{bmatrix} 1 \\ \eta \\ \mu \\ \eta\mu \end{bmatrix} \begin{bmatrix} N_2 & N_3 & N_7 & N_6 \end{bmatrix} J_{face2376} d\eta d\mu$$

$$J_{face2376} = \frac{bc}{4}$$

Detail derivation for force vector obtained by Mathematica can be found in the Appendix, it is also shown below:

$$\mathbf{f}_x = \begin{bmatrix} f_{x2} \\ f_{x3} \\ f_{x7} \\ f_{x6} \end{bmatrix} = \begin{bmatrix} \frac{1}{9}bc(9\beta_1 - 3\beta_2 - 3\beta_3 + \beta_4) \\ \frac{1}{9}bc(9\beta_1 + 3\beta_2 - 3\beta_3 - \beta_4) \\ \frac{1}{9}bc(9\beta_1 + 3\beta_2 + 3\beta_3 + \beta_4) \\ \frac{1}{9}bc(9\beta_1 - 3\beta_2 + 3\beta_3 - \beta_4) \end{bmatrix} \quad (2)$$

Similarly,

$$t_y = \sigma_{12} = \beta_{13} + \beta_{14}\mu = \begin{bmatrix} \beta_{13} & \beta_{14} \end{bmatrix} \begin{bmatrix} 1 \\ \mu \end{bmatrix}$$

$$t_z = \sigma_{13} = \beta_{17} + \beta_{18}\eta = \begin{bmatrix} \beta_{17} & \beta_{18} \end{bmatrix} \begin{bmatrix} 1 \\ \eta \end{bmatrix}$$

Forces in y and z directions for each node (2, 3, 7, and 6 in this case) can be obtained in the same manner.

$$\mathbf{f}_y = \begin{bmatrix} f_{y2} \\ f_{y3} \\ f_{y7} \\ f_{y6} \end{bmatrix} = \begin{bmatrix} bc(\beta_{13} - \frac{\beta_{14}}{3}) \\ bc(\beta_{13} - \frac{\beta_{14}}{3}) \\ bc(\beta_{13} + \frac{\beta_{14}}{3}) \\ bc(\beta_{13} + \frac{\beta_{14}}{3}) \end{bmatrix} \quad (3)$$

$$\mathbf{f}_z = \begin{bmatrix} f_{z2} \\ f_{z3} \\ f_{z7} \\ f_{z6} \end{bmatrix} = \begin{bmatrix} bc(\beta_{17} - \frac{\beta_{18}}{3}) \\ bc(\beta_{17} + \frac{\beta_{18}}{3}) \\ bc(\beta_{17} + \frac{\beta_{18}}{3}) \\ bc(\beta_{17} - \frac{\beta_{18}}{3}) \end{bmatrix} \quad (4)$$

The assembly process of *equation (2)*, *equation (3)*, and *equation (4)* has to be done to get the total force vector for four conner nodes 2, 3, 7, and 6 in *face 2376*:

$$\mathbf{f}_{face2376} = \begin{bmatrix} f_{x2} \\ f_{y2} \\ f_{z2} \\ f_{x3} \\ f_{y3} \\ f_{z3} \\ f_{x7} \\ f_{y7} \\ f_{z7} \\ f_{x6} \\ f_{y6} \\ f_{z6} \end{bmatrix} = \begin{bmatrix} \frac{1}{9}bc(9\beta_1 - 3\beta_2 - 3\beta_3 + \beta_4) \\ bc(\beta_{13} - \frac{\beta_{14}}{3}) \\ bc(\beta_{17} - \frac{\beta_{18}}{3}) \\ \frac{1}{9}bc(9\beta_1 + 3\beta_2 - 3\beta_3 - \beta_4) \\ bc(\beta_{13} - \frac{\beta_{14}}{3}) \\ bc(\beta_{17} + \frac{\beta_{18}}{3}) \\ \frac{1}{9}bc(9\beta_1 + 3\beta_2 + 3\beta_3 + \beta_4) \\ bc(\beta_{13} + \frac{\beta_{14}}{3}) \\ bc(\beta_{17} + \frac{\beta_{18}}{3}) \\ \frac{1}{9}bc(9\beta_1 - 3\beta_2 + 3\beta_3 - \beta_4) \\ bc(\beta_{13} + \frac{\beta_{14}}{3}) \\ bc(\beta_{17} - \frac{\beta_{18}}{3}) \end{bmatrix} \qquad (5)$$

Note that nodal force vector for *face 1485*, $\mathbf{f}_{face1485}$ will have the opposite sign of $\mathbf{f}_{face2376}$, and nodal force vectors for the other two pairs of faces can be done in the same method.

To this end, matrix **A**-the equilibrium matrix-is obtained after assembly of 6 nodal force vectors as done in *equation (5)*, in the FEM literature, **A** is also known as the leverage matrix.

$$f = A/3 \qquad (6)$$

$|A_z|$ matrix can be seen as following:
$24 \to 18$

$$A = \begin{bmatrix}
-\frac{bc}{4} & \frac{bc}{12} & \frac{bc}{12} & -\frac{bc}{36} & 0 & 0 & 0 & 0 & 0 & 0 & 0 & 0 & -\frac{ac}{4} & \frac{ac}{12} & 0 & 0 & -\frac{ab}{4} & \frac{ab}{12} \\
0 & 0 & 0 & 0 & -\frac{ac}{4} & \frac{ac}{12} & \frac{ac}{12} & -\frac{ac}{36} & 0 & 0 & 0 & 0 & -\frac{bc}{4} & \frac{bc}{12} & -\frac{ab}{4} & \frac{ab}{12} & 0 & 0 \\
0 & 0 & 0 & 0 & 0 & 0 & 0 & 0 & -\frac{ab}{4} & \frac{ab}{12} & \frac{ab}{12} & -\frac{ab}{36} & 0 & 0 & -\frac{ac}{4} & \frac{ac}{12} & -\frac{bc}{4} & \frac{bc}{12} \\
\frac{bc}{4} & -\frac{bc}{12} & -\frac{bc}{12} & \frac{bc}{36} & 0 & 0 & 0 & 0 & 0 & 0 & 0 & 0 & -\frac{ac}{4} & \frac{ac}{12} & 0 & 0 & -\frac{ab}{4} & -\frac{ab}{12} \\
0 & 0 & 0 & 0 & -\frac{ac}{4} & -\frac{ac}{12} & \frac{ac}{12} & \frac{ac}{36} & 0 & 0 & 0 & 0 & \frac{bc}{4} & -\frac{bc}{12} & -\frac{ab}{4} & -\frac{ab}{12} & 0 & 0 \\
0 & 0 & 0 & 0 & 0 & 0 & 0 & 0 & -\frac{ab}{4} & -\frac{ab}{12} & \frac{ab}{12} & \frac{ab}{36} & 0 & 0 & -\frac{ac}{4} & -\frac{ac}{12} & \frac{bc}{4} & -\frac{bc}{12} \\
\frac{bc}{4} & \frac{bc}{12} & -\frac{bc}{12} & -\frac{bc}{36} & 0 & 0 & 0 & 0 & 0 & 0 & 0 & 0 & \frac{ac}{4} & -\frac{ac}{12} & 0 & 0 & -\frac{ab}{4} & \frac{ab}{12} \\
0 & 0 & 0 & 0 & \frac{ac}{4} & \frac{ac}{12} & -\frac{ac}{12} & -\frac{ac}{36} & 0 & 0 & 0 & 0 & \frac{bc}{4} & -\frac{bc}{12} & -\frac{ab}{4} & \frac{ab}{12} & 0 & 0 \\
0 & 0 & 0 & 0 & 0 & 0 & 0 & 0 & -\frac{ab}{4} & -\frac{ab}{12} & -\frac{ab}{12} & -\frac{ab}{36} & 0 & 0 & \frac{ac}{4} & \frac{ac}{12} & \frac{bc}{4} & \frac{bc}{12} \\
-\frac{bc}{4} & -\frac{bc}{12} & \frac{bc}{12} & \frac{bc}{36} & 0 & 0 & 0 & 0 & 0 & 0 & 0 & 0 & \frac{ac}{4} & -\frac{ac}{12} & 0 & 0 & -\frac{ab}{4} & \frac{ab}{12} \\
0 & 0 & 0 & 0 & \frac{ac}{4} & -\frac{ac}{12} & -\frac{ac}{12} & \frac{ac}{36} & 0 & 0 & 0 & 0 & -\frac{bc}{4} & \frac{bc}{12} & -\frac{ab}{4} & \frac{ab}{12} & 0 & 0 \\
0 & 0 & 0 & 0 & 0 & 0 & 0 & 0 & -\frac{ab}{4} & \frac{ab}{12} & -\frac{ab}{12} & \frac{ab}{36} & 0 & 0 & \frac{ac}{4} & \frac{ac}{12} & -\frac{bc}{4} & -\frac{bc}{12} \\
-\frac{bc}{4} & \frac{bc}{12} & -\frac{bc}{12} & \frac{bc}{36} & 0 & 0 & 0 & 0 & 0 & 0 & 0 & 0 & -\frac{ac}{4} & -\frac{ac}{12} & 0 & 0 & \frac{ab}{4} & -\frac{ab}{12} \\
0 & 0 & 0 & 0 & -\frac{ac}{4} & \frac{ac}{12} & -\frac{ac}{12} & \frac{ac}{36} & 0 & 0 & 0 & 0 & -\frac{bc}{4} & -\frac{bc}{12} & \frac{ab}{4} & -\frac{ab}{12} & 0 & 0 \\
0 & 0 & 0 & 0 & 0 & 0 & 0 & 0 & \frac{ab}{4} & -\frac{ab}{12} & -\frac{ab}{12} & \frac{ab}{36} & 0 & 0 & -\frac{ac}{4} & -\frac{ac}{12} & -\frac{bc}{4} & \frac{bc}{12} \\
\frac{bc}{4} & -\frac{bc}{12} & \frac{bc}{12} & -\frac{bc}{36} & 0 & 0 & 0 & 0 & 0 & 0 & 0 & 0 & \frac{ac}{4} & \frac{ac}{12} & 0 & 0 & \frac{ab}{4} & \frac{ab}{12} \\
0 & 0 & 0 & 0 & \frac{ac}{4} & -\frac{ac}{12} & \frac{ac}{12} & -\frac{ac}{36} & 0 & 0 & 0 & 0 & \frac{bc}{4} & \frac{bc}{12} & \frac{ab}{4} & \frac{ab}{12} & 0 & 0 \\
0 & 0 & 0 & 0 & 0 & 0 & 0 & 0 & \frac{ab}{4} & \frac{ab}{12} & \frac{ab}{12} & \frac{ab}{36} & 0 & 0 & \frac{ac}{4} & \frac{ac}{12} & \frac{bc}{4} & \frac{bc}{12} \\
-\frac{bc}{4} & -\frac{bc}{12} & -\frac{bc}{12} & -\frac{bc}{36} & 0 & 0 & 0 & 0 & 0 & 0 & 0 & 0 & \frac{ac}{4} & \frac{ac}{12} & 0 & 0 & \frac{ab}{4} & -\frac{ab}{12} \\
0 & 0 & 0 & 0 & \frac{ac}{4} & -\frac{ac}{12} & \frac{ac}{12} & -\frac{ac}{36} & 0 & 0 & 0 & 0 & -\frac{bc}{4} & -\frac{bc}{12} & \frac{ab}{4} & -\frac{ab}{12} & 0 & 0 \\
0 & 0 & 0 & 0 & 0 & 0 & 0 & 0 & \frac{ab}{4} & -\frac{ab}{12} & \frac{ab}{12} & -\frac{ab}{36} & 0 & 0 & \frac{ac}{4} & \frac{ac}{12} & -\frac{bc}{4} & -\frac{bc}{12}
\end{bmatrix}$$

$$(7)$$

# 3 Generalized stiffness matrix

Integrate the complementary energy density over the body of the element, a flexibility $\mathbf{F}_{/3}$ matrix is then found as following:

$$U = \int_v \frac{1}{2}\boldsymbol{\sigma}^T \mathbf{C}\boldsymbol{\sigma} dV = \frac{1}{2}\boldsymbol{\beta}^T \mathbf{F}_\beta \boldsymbol{\beta},$$

then

$$\mathbf{F}_\beta = \int_v \boldsymbol{\beta}^T \boldsymbol{\sigma}^T \mathbf{C}\boldsymbol{\sigma}\boldsymbol{\beta}^{-1} dV \tag{8}$$

With *equation (1)*: $\mathbf{N} = \boldsymbol{\sigma}\boldsymbol{\beta}^{-1}$ previously defined from the assumed stress field, the $\mathbf{F}_\beta$ in *equation (8)* can be calculated as:

$$\mathbf{F}_\beta = \int_v \mathbf{N}^T \mathbf{C} \mathbf{N} dV \tag{9}$$

Its inverse will be the generalized stiffness matrix $\mathbf{S}_\beta = \mathbf{F}_\beta^{-1}$:

$$\underbrace{\mathbf{S}_\beta}_{18\times18} = \begin{bmatrix} \underbrace{\mathbf{S}_{\beta,11}}_{18\times3} & \underbrace{\mathbf{S}_{\beta,12}}_{18\times3} & \underbrace{\mathbf{S}_{\beta,13}}_{18\times3} & \underbrace{\mathbf{S}_{\beta,14}}_{18\times3} & \underbrace{\mathbf{S}_{\beta,15}}_{18\times6} \end{bmatrix} \tag{10}$$

An entire $\mathbf{S}_{/3}$ matrix could be seen in the enclosed Appendix, while its first three columns are shown below:

$$S_{\beta,11} = \begin{bmatrix} -\frac{a^2b^2c^2}{Em^2} - \frac{a^2b^2c^2\nu^2}{Em^2} - \frac{2a^3b^3\nu^3c^3}{Em^3} - \frac{3a^3b^3\nu^2c^3}{Em^3} + \frac{a^3b^3c^3}{Em^3} & 0 & 0 \\ 0 & \frac{abc}{3Em\left(\frac{a^2b^2c^2}{9Em^2} - \frac{a^2b^2c^2\nu^2}{9Em^2}\right)} & 0 \\ 0 & 0 & \frac{abc}{3Em\left(\frac{a^2b^2c^2}{9Em^2} - \frac{a^2b^2c^2\nu^2}{9Em^2}\right)} \\ 0 & 0 & 0 \\ \frac{a^2b^2\nu^2c^2}{Em^2} + \frac{a^2b^2\nu c^2}{Em^2} - \frac{2a^3b^3\nu^3c^3}{Em^3} - \frac{3a^3b^3\nu^2c^3}{Em^3} + \frac{a^3b^3c^3}{Em^3} & 0 & 0 \\ 0 & 0 & 0 \\ 0 & 0 & \frac{abc\nu}{3Em\left(\frac{a^2b^2c^2}{9Em^2} - \frac{a^2b^2c^2\nu^2}{9Em^2}\right)} \\ 0 & 0 & 0 \\ \frac{a^2b^2\nu^2c^2}{Em^2} + \frac{a^2b^2\nu c^2}{Em^2} - \frac{2a^3b^3\nu^3c^3}{Em^3} - \frac{3a^3b^3\nu^2c^3}{Em^3} + \frac{a^3b^3c^3}{Em^3} & 0 & 0 \\ 0 & 0 & 0 \\ 0 & \frac{abc\nu}{3Em\left(\frac{a^2b^2c^2}{9Em^2} - \frac{a^2b^2c^2\nu^2}{9Em^2}\right)} & 0 \\ 0 & 0 & 0 \\ 0 & 0 & 0 \\ 0 & 0 & 0 \\ 0 & 0 & 0 \\ 0 & 0 & 0 \\ 0 & 0 & 0 \\ 0 & 0 & 0 \end{bmatrix}$$

# 4 Physical stiffness matrix

The physical stiffness $\mathbf{K}_{cr}$ relates $\mathbf{f} = \mathbf{K}_{cr}\mathbf{u}$. From previous derivation, *equation (5)*: $\mathbf{f} = \mathbf{A}/3 = \mathbf{S}_{/3}\mathbf{A}^T\mathbf{u}$ which can give:

$$\mathbf{K}_{cr} = \mathbf{A}\mathbf{S}_{/3}\mathbf{A}^T.$$

Since $\mathbf{K}_{cr}$ is quite large as the result of the symbolic computation, it can

$$\underbrace{\mathbf{K}_\sigma}_{24\times 24} = \left[ \underbrace{\mathbf{K}_{\sigma,1}}_{24\times 2}\ \underbrace{\mathbf{K}_{\sigma,2}}_{24\times 2}\ \underbrace{\mathbf{K}_{\sigma,3}}_{24\times 2}\ \underbrace{\mathbf{K}_{\sigma,4}}_{24\times 2}\ \underbrace{\mathbf{K}_{\sigma,5}}_{24\times 2}\ \underbrace{\mathbf{K}_{\sigma,6}}_{24\times 2}\ \underbrace{\mathbf{K}_{\sigma,7}}_{24\times 2}\ \underbrace{\mathbf{K}_{\sigma,8}}_{24\times 2}\ \underbrace{\mathbf{K}_{\sigma,9}}_{24\times 2}\ \underbrace{\mathbf{K}_{\sigma,10}}_{24\times 2}\ \underbrace{\mathbf{K}_{\sigma,11}}_{24\times 2}\ \underbrace{\mathbf{K}_{\sigma,12}}_{24\times 2} \right]$$

and the first two columns of $\mathbf{K}_{cr}$ are as following:

$$\mathbf{K}_{\sigma,1} = \begin{bmatrix} \dots & \dots \end{bmatrix}$$

As seen in the provided Appendix, the rank of $\mathbf{K}_{cr}$ is 18 which is correct since there are 24 degrees of freedom for 8-node brick element and 6 rigid body motions.

## 5  Higher order stiffness matrix

The stiffness matrix $\mathbf{K}_{cr}$ can be decomposed as:

$$\underbrace{\mathbf{K}_\sigma}_{24\times 24} = V^{-1} \underbrace{\mathbf{L}}_{24\times 6} \underbrace{\mathbf{E}}_{6\times 6} \underbrace{\mathbf{L}^T}_{6\times 24} + V \underbrace{\mathbf{H}_h^T}_{24\times 12} \underbrace{\mathbf{W}^T}_{12\times 12} \underbrace{\mathbf{R}}_{12\times 12} \underbrace{\mathbf{W}}_{12\times 12} \underbrace{\mathbf{H}_h}_{12\times 24} \qquad (11)$$

$$= \mathbf{K}_b + \mathbf{K}_h \qquad (12)$$

where V=abc is the element volume, $\mathbf{H}_h$, and $\mathbf{L}$ are the same for stiffness matrices based on three different methods, i.e. assumed stress method, assumed strain method, and assumed displacement method. The force lumping $\mathbf{L}$ matrix is as following:

$$\mathbf{L} = \begin{bmatrix}
-\frac{bc}{4} & 0 & 0 & -\frac{ac}{4} & 0 & -\frac{ab}{4} \\
0 & -\frac{ac}{4} & 0 & -\frac{bc}{4} & -\frac{ab}{4} & 0 \\
0 & 0 & -\frac{ab}{4} & 0 & -\frac{ac}{4} & -\frac{bc}{4} \\
\frac{bc}{4} & 0 & 0 & -\frac{ac}{4} & 0 & -\frac{ab}{4} \\
0 & -\frac{ac}{4} & 0 & \frac{bc}{4} & -\frac{ab}{4} & 0 \\
0 & 0 & -\frac{ab}{4} & 0 & -\frac{ac}{4} & \frac{bc}{4} \\
\frac{bc}{4} & 0 & 0 & \frac{ac}{4} & 0 & -\frac{ab}{4} \\
0 & \frac{ac}{4} & 0 & \frac{bc}{4} & -\frac{ab}{4} & 0 \\
0 & 0 & -\frac{ab}{4} & 0 & \frac{ac}{4} & \frac{bc}{4} \\
-\frac{bc}{4} & 0 & 0 & \frac{ac}{4} & 0 & -\frac{ab}{4} \\
0 & \frac{ac}{4} & 0 & -\frac{bc}{4} & -\frac{ab}{4} & 0 \\
0 & 0 & -\frac{ab}{4} & 0 & \frac{ac}{4} & -\frac{bc}{4} \\
-\frac{bc}{4} & 0 & 0 & -\frac{ac}{4} & 0 & \frac{ab}{4} \\
0 & -\frac{ac}{4} & 0 & -\frac{bc}{4} & \frac{ab}{4} & 0 \\
0 & 0 & \frac{ab}{4} & 0 & -\frac{ac}{4} & -\frac{bc}{4} \\
\frac{bc}{4} & 0 & 0 & -\frac{ac}{4} & 0 & \frac{ab}{4} \\
0 & -\frac{ac}{4} & 0 & \frac{bc}{4} & \frac{ab}{4} & 0 \\
0 & 0 & \frac{ab}{4} & 0 & -\frac{ac}{4} & \frac{bc}{4} \\
\frac{bc}{4} & 0 & 0 & \frac{ac}{4} & 0 & \frac{ab}{4} \\
0 & \frac{ac}{4} & 0 & \frac{bc}{4} & \frac{ab}{4} & 0 \\
0 & 0 & \frac{ab}{4} & 0 & \frac{ac}{4} & \frac{bc}{4} \\
-\frac{bc}{4} & 0 & 0 & \frac{ac}{4} & 0 & \frac{ab}{4} \\
0 & \frac{ac}{4} & 0 & -\frac{bc}{4} & \frac{ab}{4} & 0 \\
0 & 0 & \frac{ab}{4} & 0 & \frac{ac}{4} & -\frac{bc}{4}
\end{bmatrix}$$

For isotropic materials, elasticity $\mathbf{E}$ is shown below:

$$\mathbf{E} = \begin{bmatrix} \frac{Em(1-\nu)}{(1-2\nu)(\nu+1)} & \frac{Em\nu}{(1-2\nu)(\nu+1)} & \frac{Em\nu}{(1-2\nu)(\nu+1)} & 0 & 0 & 0 \\ \frac{Em\nu}{(1-2\nu)(\nu+1)} & \frac{Em(1-\nu)}{(1-2\nu)(\nu+1)} & \frac{Em\nu}{(1-2\nu)(\nu+1)} & 0 & 0 & 0 \\ \frac{Em\nu}{(1-2\nu)(\nu+1)} & \frac{Em\nu}{(1-2\nu)(\nu+1)} & \frac{Em(1-\nu)}{(1-2\nu)(\nu+1)} & 0 & 0 & 0 \\ 0 & 0 & 0 & \frac{Em\left(\frac{1}{2}-\nu\right)}{(1-2\nu)(\nu+1)} & 0 & 0 \\ 0 & 0 & 0 & 0 & \frac{Em\left(\frac{1}{2}-\nu\right)}{(1-2\nu)(\nu+1)} & 0 \\ 0 & 0 & 0 & 0 & 0 & \frac{Em\left(\frac{1}{2}-\nu\right)}{(1-2\nu)(\nu+1)} \end{bmatrix}$$

It is important to point out in advance here that the geometric projector $\mathbf{H}_h$ matrix and a latter discussed $\mathbf{G}_{rc}$ matrix should satisfy the orthogonal check. This requirement shall be covered subsequently, for now, as seen in the second part of the summation for $\mathbf{K}_\sigma$, $\mathbf{H}_h$ matrix is then introduced as:

$$\mathbf{H}_h^T = \begin{bmatrix}
-\frac{ab}{4} & -\frac{ac}{4} & 0 & 0 & 0 & 0 & \frac{bc}{4} & 0 & 0 & -\frac{1}{8}abc & 0 & 0 \\
0 & 0 & -\frac{ab}{4} & -\frac{bc}{4} & 0 & 0 & 0 & \frac{ac}{4} & 0 & 0 & -\frac{1}{8}abc & 0 \\
0 & 0 & 0 & 0 & -\frac{ac}{4} & -\frac{bc}{4} & 0 & 0 & \frac{ab}{4} & 0 & 0 & -\frac{1}{8}abc \\
\frac{ab}{4} & \frac{ac}{4} & 0 & 0 & 0 & 0 & \frac{bc}{4} & 0 & 0 & \frac{abc}{8} & 0 & 0 \\
0 & 0 & \frac{ab}{4} & -\frac{bc}{4} & 0 & 0 & 0 & -\frac{ac}{4} & 0 & 0 & \frac{abc}{8} & 0 \\
0 & 0 & 0 & 0 & \frac{ac}{4} & -\frac{bc}{4} & 0 & 0 & -\frac{ab}{4} & 0 & 0 & \frac{abc}{8} \\
-\frac{ab}{4} & \frac{ac}{4} & 0 & 0 & 0 & 0 & -\frac{bc}{4} & 0 & 0 & -\frac{1}{8}abc & 0 & 0 \\
0 & 0 & -\frac{ab}{4} & \frac{bc}{4} & 0 & 0 & 0 & -\frac{ac}{4} & 0 & 0 & -\frac{1}{8}abc & 0 \\
0 & 0 & 0 & 0 & \frac{ac}{4} & \frac{bc}{4} & 0 & 0 & \frac{ab}{4} & 0 & 0 & -\frac{1}{8}abc \\
\frac{ab}{4} & -\frac{ac}{4} & 0 & 0 & 0 & 0 & -\frac{bc}{4} & 0 & 0 & \frac{abc}{8} & 0 & 0 \\
0 & 0 & \frac{ab}{4} & \frac{bc}{4} & 0 & 0 & 0 & \frac{ac}{4} & 0 & 0 & \frac{abc}{8} & 0 \\
0 & 0 & 0 & 0 & -\frac{ac}{4} & \frac{bc}{4} & 0 & 0 & -\frac{ab}{4} & 0 & 0 & \frac{abc}{8} \\
-\frac{ab}{4} & \frac{ac}{4} & 0 & 0 & 0 & 0 & -\frac{bc}{4} & 0 & 0 & \frac{abc}{8} & 0 & 0 \\
0 & 0 & -\frac{ab}{4} & \frac{bc}{4} & 0 & 0 & 0 & -\frac{ac}{4} & 0 & 0 & \frac{abc}{8} & 0 \\
0 & 0 & 0 & 0 & \frac{ac}{4} & \frac{bc}{4} & 0 & 0 & \frac{ab}{4} & 0 & 0 & \frac{abc}{8} \\
\frac{ab}{4} & -\frac{ac}{4} & 0 & 0 & 0 & 0 & -\frac{bc}{4} & 0 & 0 & -\frac{1}{8}abc & 0 & 0 \\
0 & 0 & \frac{ab}{4} & \frac{bc}{4} & 0 & 0 & 0 & \frac{ac}{4} & 0 & 0 & -\frac{1}{8}abc & 0 \\
0 & 0 & 0 & 0 & -\frac{ac}{4} & \frac{bc}{4} & 0 & 0 & -\frac{ab}{4} & 0 & 0 & -\frac{1}{8}abc \\
-\frac{ab}{4} & -\frac{ac}{4} & 0 & 0 & 0 & 0 & \frac{bc}{4} & 0 & 0 & \frac{abc}{8} & 0 & 0 \\
0 & 0 & -\frac{ab}{4} & -\frac{bc}{4} & 0 & 0 & 0 & \frac{ac}{4} & 0 & 0 & \frac{abc}{8} & 0 \\
0 & 0 & 0 & 0 & -\frac{ac}{4} & -\frac{bc}{4} & 0 & 0 & \frac{ab}{4} & 0 & 0 & \frac{abc}{8} \\
\frac{ab}{4} & \frac{ac}{4} & 0 & 0 & 0 & 0 & \frac{bc}{4} & 0 & 0 & -\frac{1}{8}abc & 0 & 0 \\
0 & 0 & \frac{ab}{4} & -\frac{bc}{4} & 0 & 0 & 0 & -\frac{ac}{4} & 0 & 0 & -\frac{1}{8}abc & 0 \\
0 & 0 & 0 & 0 & \frac{ac}{4} & -\frac{bc}{4} & 0 & 0 & -\frac{ab}{4} & 0 & 0 & -\frac{1}{8}abc
\end{bmatrix}$$

According to Felippa's notation in the study in 2006, **W** is a higher order mode weighting matrix. In particular, for a rectangular panel, **W** is diagonal and formulation independent (Felippa, 2006). Shown below is **W** matrix:

$$\mathbf{W} = \begin{bmatrix} \frac{1}{ab} & 0 & 0 & 0 & 0 & 0 & 0 & 0 & 0 & 0 & 0 & 0 \\ 0 & \frac{1}{ac} & 0 & 0 & 0 & 0 & 0 & 0 & 0 & 0 & 0 & 0 \\ 0 & 0 & \frac{1}{ab} & 0 & 0 & 0 & 0 & 0 & 0 & 0 & 0 & 0 \\ 0 & 0 & 0 & \frac{1}{bc} & 0 & 0 & 0 & 0 & 0 & 0 & 0 & 0 \\ 0 & 0 & 0 & 0 & \frac{1}{ac} & 0 & 0 & 0 & 0 & 0 & 0 & 0 \\ 0 & 0 & 0 & 0 & 0 & 1 & 0 & 0 & 0 & 0 & 0 & 0 \\ 0 & 0 & 0 & 0 & 0 & 0 & \frac{1}{bc} & 0 & 0 & 0 & 0 & 0 \\ 0 & 0 & 0 & 0 & 0 & 0 & 0 & \frac{1}{a} & 0 & 0 & 0 & 0 \\ 0 & 0 & 0 & 0 & 0 & 0 & 0 & 0 & \frac{1}{ab} & 0 & 0 & 0 \\ 0 & 0 & 0 & 0 & 0 & 0 & 0 & 0 & 0 & \frac{1}{abc} & 0 & 0 \\ 0 & 0 & 0 & 0 & 0 & 0 & 0 & 0 & 0 & 0 & \frac{1}{abc} & 0 \\ 0 & 0 & 0 & 0 & 0 & 0 & 0 & 0 & 0 & 0 & 0 & \frac{1}{abc} \end{bmatrix}$$

The basic stiffness matrix $\mathbf{K}_b$ that takes care of consistency and mixability can be calculated and whose rank should be 6:

$$\mathbf{K}_b = V^{-1} \mathbf{LEL}^T$$

Here, it is interesting to be based upon the aforementioned *equation (12)* which allows one to find the high order stiffness matrix. As the result, the higher order stiffness matrix $\mathbf{K}_h$ that provides correct rank and accuracy can be calculated:

$$\mathbf{K}_h = \mathbf{K}_{cr} - \mathbf{K}_b$$

$$\{\mathbf{K}_h\}_{24 \to 24} = \begin{bmatrix} \mathbf{K}_{h,1} & \mathbf{K}_{h,2} & \mathbf{K}_{h,3} & \ldots & \mathbf{K}_{h,22} & \mathbf{K}_{h,23} & \mathbf{K}_{h,24} \\ \{z\}_{24 \to 1} & \{z\}_{24 \to 1} & \{z\}_{24 \to 1} & & \{z\}_{24 \to 1} & \{z\}_{24 \to 1} & \{z\}_{24 \to 1} \end{bmatrix}$$

The following solely displays the first column of $\mathbf{K}_h$:

The rank of $\mathbf{K}_h$ is 12 which is correct since $\text{rank}(\mathbf{K}_{cr}) - \text{rank}(\mathbf{K}_b) = 18 - 6 = 12$.

When the higher order matrix $\mathbf{K}_h$ provides the correct rank, it means that $\mathbf{K}_h$ satisfies the Stability requirement. More interesting, $\mathbf{K}_h$ is orthogonal to rigid body motions and constant strain states. This claim leads to an introduction of the basic mode matrix $\mathbf{G}_{rc}$ which its columns should span the rigid body modes and constant strain states evaluated at the nodes. $\mathbf{G}_{rc}$ is shown below:

$$\begin{bmatrix} 1 & 0 & 0 & \frac{139}{1600} - \frac{b}{2} & 0 & \frac{c}{2} - \frac{351}{3200} \end{bmatrix}$$

With given $\mathbf{H}_h, \mathbf{G}_{rc}$, and Mathematica computation, it is easy to see that $\mathbf{H}_h \mathbf{G}_{rc} = 0$.

# 6  Bending Test

Knowing that there are 6 bending modes in the HO Stiffness Matrix, the bending test involves checking the in-plane bending for each side of the brick element. This leads to a comparison in the form of energy ratio with respect to each plane. This methodology is adopted from Felippa's study in 2006 (Felippa, 2006). There are 6 planes in a cube, each plane is subjected to (2) bending modes correspond to the (2) axis of the plane. For each plane, the exact solution is as for bending in beam theory. The stress in x plane is the associated stress field $\sigma_{xx} = \frac{-M_x y}{I_b}$, $\sigma_{yy} = \sigma_{xy} = 0$. With $I_b = \frac{hb^3}{12}$. For the y bending test, the beam with a cross section has height **a** and thickness **h**. Associated stress field will be $\sigma_{yy} = \frac{M_y x}{I_a}$, $\sigma_{xx} = \sigma_{xy} = 0$. With $I_a = \frac{ha^3}{12}$

With the FEM discretization, $U_x = \frac{6 a C_{11} M_x^2}{b^3 h}$, while $U_y = \frac{6 b C_{22} M_y^2}{a^3 h}$

Listed below is a $\mathbf{u}_b$ vector which contains nodal displacements due to bending:

The strain energies absorbed by the panel element under these applied nodal displacement are:

$$U_x^{panel} = \frac{\mathbf{u}^T \mathbf{Ku}_{bx}}{2} \quad \text{and} \quad U_y^{panel} = \frac{\mathbf{u}_{by}^T \mathbf{Ku}_{by}}{2}$$

The energy ratios are defined as:

$$r_x = \frac{U_x^{panel}}{U_{m\,x}^{bea}} \quad \text{and} \quad r_y = \frac{U_y^{panel}}{U_{m\,y}^{bea}}$$

when $r_x = r_y = 1$, one should expect an exact answer under bending, in other words, the panel is in plane x-bending and y-bending exact. When $r_x$ or $r_y > 1$, the panel is oversti⌋ in x bending or y bending, respectively. Contradictorily, when $r_x$ or $r_y < 1$, the panel is overflexible in x bending or y bending, respectively. Results from the bending test will be listed in the following section, where energy is calculated with $\mathbf{K}_{cr}$ as well as with $\mathbf{K}_h$.

### Energy ratio with $\mathbf{K}_{cr}$

Energy ratio with $\mathbf{K}_{cr}$ will be:

$$r_{xy} = \frac{a^2(\ -1) - 32c^2}{32c^2(\ ^2 - 1)} \tag{13}$$

The following figure shows the path of the energy ratio when the height of the panel is fixed, and the length of the panel is set from 1 to 10. The x axis is defined as "Aspect ratio" which is the ratio of the panel length over the panel height, while the y axis is called "Energy ratio" as discussed previously.

Figure 2. Energy ratio with **K**$_{o\text{-}}$

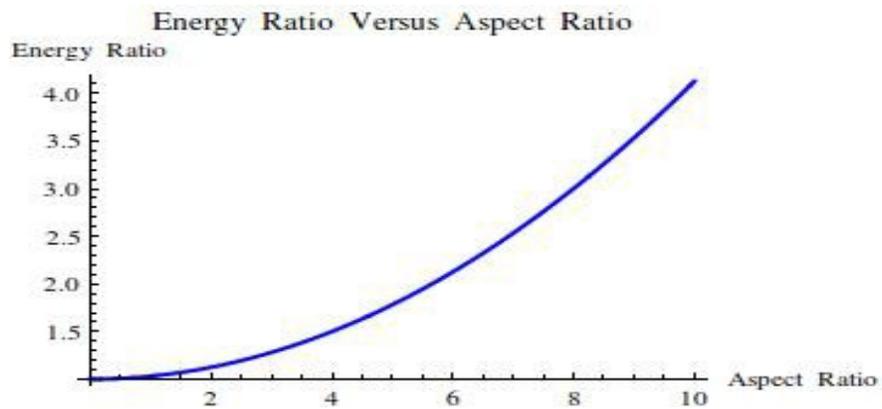

As seen in these figures above, when "Aspect ratio" ranges from 0 to 2, the element performs ideally as the result of the energy ratio runs from 1 to 1.5.

For simplification purpose, Poisson ratio was set to zero in *equation (13)*, energy ratio thus equals:

$$r_{xyo} = \frac{a^2}{32c^2} + 1$$

Figure 3. Energy ratio with **K**$_{o\text{-}}$ with zeroed Poisson ratio

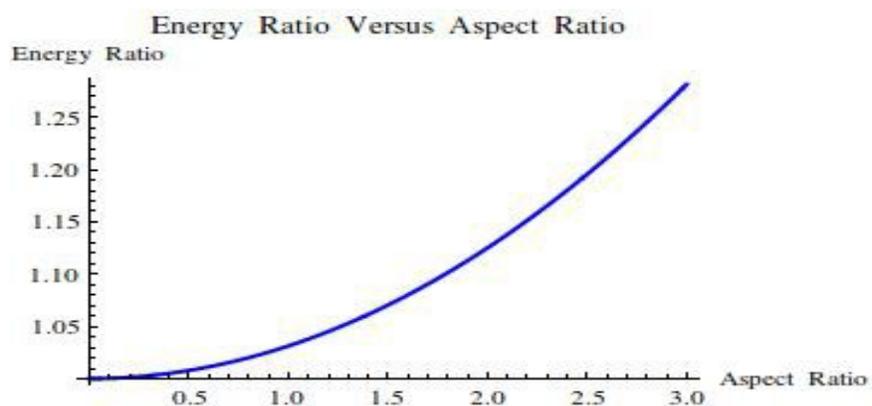

## Energy ratio with K*h*

As aforementioned, **K**$_h$ will be substituted into the energy formula, thus, seen below is the relationship of energy ratio with **K**$_h$:

$$r_{xy} = \frac{a^2(\nu - 1) - 32c^2}{32c^2(\nu^2 - 1)} \tag{14}$$

Figure 4. Energy ratio with **K**$_h$ and zero Poisson ratio

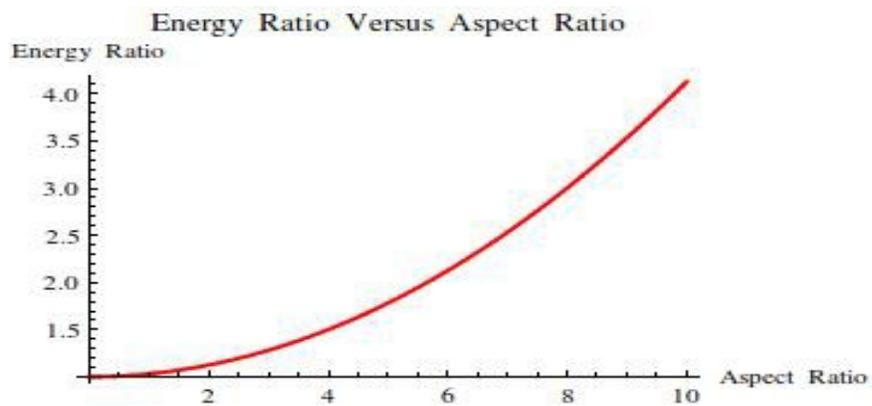

And when set Poisson ratio to zero in *equation (14)*, one should also obtain:

$$r_{xyo} = \frac{a^2}{32c^2} + 1$$

Figure 5. Energy ratio with **K**$_o$- and zero Poisson ratio

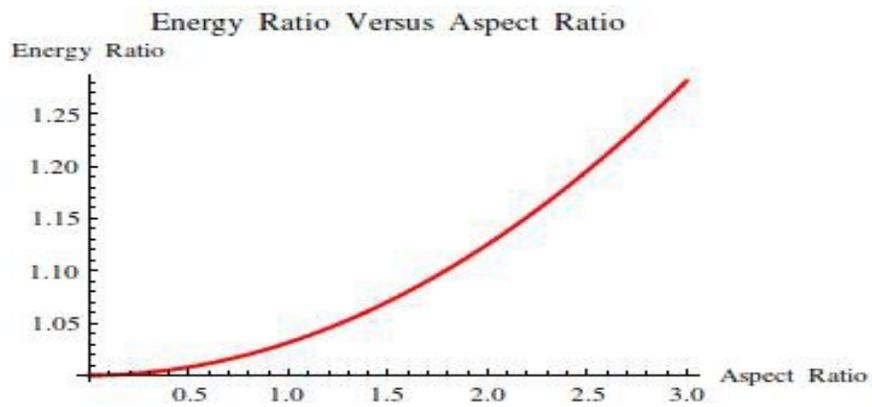

# 7 Conclusion

This study covered a development of a high-order stiffness matrix for a hexahedron brick element. The resulting benchmark test used bending and defined unity energy ratio and aspect ratio led to a methodology for an optimal element. Symbolic calculations were carried out by Mathematica which fortunately made the computation processes much more efficient throughout the complicated journeys of all variables. Numerical derivations were demonstrated in the enclosed Appendix. Finally, the contribution in terms of new results include the following:

(1) To the authors understanding, this is by far the first template for a three dimensional solid element.

(2) The study draws out a framework that described a methodology for finding an optimal stiffness matrix in FE element.

(3) The study has led to one useful extension which is an ongoing research focuses on finding a computation for optimal free parameters using mathematical programming methods, by minimizing squared deviations of energy ratios from unity.

Future study should include the effect of non-parallel geometry in a hexahedron element as well as the development of templates for 3D solid-shell elements. Additionally, a template verification using the Assumes Natural Strain method could be done for this particular element.